\documentclass[11pt]{article}
\usepackage{amsmath}
\usepackage{amssymb}
\usepackage{amsfonts} 
\usepackage{amstext} 
\usepackage{amsthm} 
\usepackage{euscript} 

\newcommand{\Z}{\ensuremath{\mathbb{Z}}}
\newcommand{\N}{\ensuremath{\mathbb{N}}}

\newcommand{\R}{\ensuremath{\mathbb{R}}} 
\newcommand{\C}{\ensuremath{\mathbb{C}}}

\newcommand{\W}{\ensuremath{\Omega}} 
\newcommand{\w}{\ensuremath{\omega}} 
\newcommand{\wc}{\ensuremath{\omega_{\text{can}}}} 
\newcommand{\lac}{\ensuremath{\lambda_{\text{can}}}} 
\newcommand{\al}{\ensuremath{\alpha}} 
\newcommand{\be}{\ensuremath{\beta}} 
\newcommand{\g}{\ensuremath{\gamma}} 
\newcommand{\de}{\ensuremath{\delta}}

\newcommand{\io}{\ensuremath{\iota}} 
\newcommand{\h}{\ensuremath{\theta}} 
\newcommand{\sg}{\ensuremath{\Sigma}} 
\newcommand{\s}{\ensuremath{\sigma}} 
\newcommand{\la}{\ensuremath{\lambda}} 
\newcommand{\Diff}{\text{Diff}} 
\newcommand{\Symp}{\text{Symp}}

\newcommand{\cE}{\ensuremath{\mathcal{E}}}

\newcommand{\cM}{\ensuremath{\mathcal{M}}} 
\newcommand{\cN}{\ensuremath{\mathcal{N}}}

\newtheorem{defn}{Definition}[section]
\newtheorem{exam}{Example}[section]
\newtheorem{lem}{Lemma}[section]
\newtheorem{cor}{Corollary}[section]

\newtheorem{theo}{Theorem}[section]

\newtheorem{rem}{Remark}[section]

\newtheorem{ques}{Question}[section]

\begin{document} 

\begin{center} 
{\Large \bf Invariants of  
Lagrangian surfaces} \\ 

\vspace{.3in} 
Mei-Lin Yau  \\  
   \vspace{.2in} 
Department of Mathematics,  Michigan State University
\\  East Lansing, MI 48824 \\ 
Email: yau@math.msu.edu
\end{center} 
\vspace{.2in} 
\begin{abstract} 
We define 
a nonnegative integer $\la(L,L_0;\phi)$ 
for a pair of diffeomorphic closed Lagrangian 
surfaces $L_0,L$ embedded in a symplectic 4-manifold 
$(M,\w)$ and a diffeomorphism $\phi \in\Diff^+(M)$ 
satisfying $\phi(L_0)=L$. We prove that if there exists 
$\phi\in\Diff^+_o(M)$ with $\phi(L_0)=L$ and 
$\la(L,L_0;\phi)=0$, then 
$L_0,L$ are symplectomorphic.  
 We also define a    
second invariant $n(L_1,L_0;[L_t])=n(L_1,L_0,[\phi_t])$ 
for a smooth isotopy $L_t=\phi_t(L_0)$ between 
two Lagrangian surfaces $L_0$ and 
$L_1$ with $\la (L_1,L_0;\phi_1)=0$, 
which serves as an obstruction of deforming 
$L_t$ to a Lagrangian isotopy with $L_0,L_1$ preserved. 
\end{abstract}

\section{Introduction} \label{intro}

One subtle question in symplectic topology is to find the 
fine line between symplectic topology and differential topology. 
For example, what are the things that can be done 
diffeomorphically 
but not symplectically? The objects to be tested on are 
embedded compact Lagrangian surfaces. In their own worlds, 
individual Lagrangian surfaces do not know the existence of 
symplectic structures, until they try to communicate with 
each other via diffeomorphisms and/or homotopies. 
Diffeomorphic Lagrangian surfaces may be surprised 
to find that they 
live in quite different neighborhoods 
dictated by symplectic structures. Even if this is not the 
case, smoothly isotopic Lagrangian surfaces may still
find that they are 
destined to meet the symplectic structure (by becoming 
non-Lagrangian) before they meet each other,  
no matter which path they choose. Such phenomena have been 
explored and studied by Fintushel-Stern \cite{FS} and 
Seidel \cite{Se} (see also \cite{EP}).

In this note we construct two invariants of 
{\em embedded} compact Lagrangian 
surfaces to address the two questions described 
above:   
\begin{itemize} 
\item When two diffeomorphic Lagrangian surfaces are 
symplectomorphic? 
\item When two smoothly isotopic Lagrangian surfaces are 
Lagrangian isotopic? 
\end{itemize} 
Here compact Lagrangian surfaces 
$L_0,L_1$ in a symplectic 4--manifold $(M,\w)$ are said to  
be smoothly isotopic if there exists 
between $L_0$ and $L_1$ a smooth homotopy consists 
of embeddings. This is equivalent to the existence of 
a smooth family $\phi_t\in \Diff^+(M)$ with $\phi_0=id$ and 
$\phi_1(L_0)=L_1$. $L_0,L_1$ are Lagrangian isotopic if 
the smooth isotopy consists of Lagrangian surfaces. 

Let $L_0,L$ be two embedded Lagrangian surface such that 
$L=\phi(L_0)$ for some $\phi\in \Diff^+(M)$. 
The first invariant $\la (L,L_0;\phi)$ 
we construct is really a generalization 
of the $\la (T)$ invariant (for Lagrangian torus $T$) 
constructed by Fintushel and Stern \cite{FS}. This invariant 
is related to the first question we mentioned above. In its 
most general form, it is really an invariant of a pair of 
diffeomorphic Lagrangian surfaces of positive genus 
(the invariant is trivial when the genus is 0) together 
with a diffeomorphism between them.  

Define 
\begin{align*}  
\Diff^+_o(M;L_0\to L)&:=\{ \phi\in \Diff^+_o(M) \mid 
\phi(L_0)=L\} \\ 
\Symp(M,\w)&:=\{ \phi\in \Diff^+(M) \mid 
\phi^*\w=\w\} 
\end{align*}  
Here $\Diff^+_o(M)$ 
denotes the connected component of $\Diff^+(M)$. 

\begin{theo} \label{lambda} 
Let $L_0,L$ be two closed Lagrangian surfaces 
embedded in a symplectic 4-manifold $(M,\w)$. 
Assume that 
$\Diff^+_o(M;L_0\to L)\neq\emptyset$. 
Then 
$\phi\in\Diff^+_o(M;L_0\to L)$  is homotopic in 
$\Diff^+_o(M;L_0\to L)$ to some 
$\psi\in \Symp(M,\w)$ iff 
$\la (L,L_0;\phi)=0$.   
\end{theo} 

\begin{rem} 
{\rm The invariant $\la(L,L_0;\cdot )$ actually assigned to 
each connected component of $\Diff^+_o(M;L_0\to L)$ a 
nonnegative integer. Theorem \ref{lambda} then says that 
this number is 0 iff the corresponding connected component 
of $\Diff^+_o(M;L_0\to L)$ contains an element 
of $\Symp(M,\w)$. 
} 
\end{rem} 

The construction of $\la(L,L_0;\phi)$ implies that 
$\la(L,L_0;\cdot )=0$ if $L_0,L$ are of genus 0, i.e., 
if $L_0,L$ are embedded Lagrangian spheres. Hence we have the 
following 

\begin{cor} 
Let $L_0,L$ be two embedded lagrangian spheres in a symplectic 
4-manifold $(M,\w)$. Suppose that $L_0,L$ are smoothly isotopic, 
then $L_0,L$ are symplectomorphic. 
\end{cor}

The second invariant $n(L_1,L_0;[L_t])$ 
(see Section \ref{second})
seems related to the generalized Dehn 
twist considered by Seidel \cite{Se}. It is an 
invariant 
of a smooth homotopy between two Lagrangian surfaces, hence 
is related to the second question above. 

The organization of this article is as follows: In Section 
\ref{first} we construct a local version of $\la$ on  
cotangent bundles $T^*L$ and then introduce the definition 
of $\la(L,L_0;\phi)$. The rest of Section \ref{first} is 
devoted to the proof of Theorem \ref{lambda}. The  
invariant $n(L_1,L_0;[L_t])$ is constructed in Section 
\ref{second}, followed by some discussions  
and final remarks. 

\vspace{.3in} 
\noindent 
{\bf \Large Acknowledgment} \\ 

We would like express our many thanks to Ronald Fintushel 
whose inspirational seminar talk on invariants of Lagrangian 
tori (see \cite{FS}) prompted our interest in this direction.  
We also also grateful to John McCarthy for very helpful 
discussions and comments.

\section{The first invariant} \label{first}

\subsection{Local definition of $\la (L,\w',\w)$} 
\label{local-1}

Let $L$ be a closed orientable surface and 
$T^*L$ denotes the cotangent bundles of $L$. 
Let $\la_{\text{can}}$ 
denote the canonical 1-form 
Let $\wc=-d\lac$ denote the canonical 
symplectic 2--form on $T^*L$. 
Let $x$ be a local coordinate on $L$ and $y\in \R^2$ be a 
local coordinate on the fiber $T^*_xL$, then $(x,y)$ are 
local coordinates of $T^*L$, and $\lac$, $\wc$ 
are given by the formulae 
\[ 
\la_{\text{can}}=ydx, \quad \wc=dx\wedge dy 
\] 
Let us fix an orientation for $T^*L$ so that $\wc^2$ becomes 
a volume form on $T^*L$. we will consider only symplectic 
2--forms $\w$ on $T^*L$ with $\w^2>0$. 
View $L$ 
as the zero section of $T^*L$. Let $\w$ be be a symplectic 
2--form on $T^*L$ with $L$ a Lagrangian surface. 
With local coordinates $(x,y)$ for $T^*L$ 
we have the following isomorphisms  
\begin{eqnarray}  
T_{(x,0)}T^*L &\cong T_xL\oplus T^*_xL  \label{iso1}\\ 
T^*_{(x,0)}T^*L &\cong T^*_xL\oplus T_xL  \label{iso2} 
\end{eqnarray}  
Being symplectic, $\w$ induces an isomorphism 
\begin{align*} 
T_{(x,y)}T^*L & \to T^*_{(x,y)}T^*L \\   
v & \to \w (v,\cdot ). 
\end{align*}  
Then by (\ref{iso1}) and (\ref{iso2}) 
$\w$ induces a bundle isomorphism 
\[ 
\W_{\w}:T_xL \to T_x^*L, \quad  
\W (v):=\w (v,\cdot ) \text{ for } v\in T_xL .  
\] 
If we fix a Riemannian metric on $L$ and let $g$ denote the 
induced Riemannian metric on $T^*L$ the cotangent bundle, 
then the splitting in (\ref{iso1}) is orthogonal. 

Let $\w'$ be another symplectic 2-form on $T^*L$ with 
$L$ Lagrangian, then 
\[ 
\W_{\w',\w}:=\W_{w'}\circ \W_{\w}^{-1}: T^*_xL \to T_x^*L 
\] 
is an orientation-preserving linear automorphism for every $x\in
L$,  hence is a section of the trivial $GL(2,\R)$-bundle over 
$L$. By fixing a trivialization of the bundle we get 
\begin{equation} \label{omega}  
\W_{\w',\w}: L\to GL(2,\R). 
\end{equation}  
Since $GL(2,\R)$ is homotopic to $S^1$, the homotopy class of 
the map (\ref{omega}) is classified by 
\[ 
[L,S^1]\cong H^1(L,\Z)\cong \Z^{2g}  
\] 
where $g$ is the genus of $L$, 
and is independent of the choice of the trivialization of 
the trivial $GL(2,\R)$-bundle. Note that the isomorphism 
$H^1(L,\Z)\cong \Z^{2g}$ depends on a choice of a basis 
for $H^1(L,\Z)\cong \Z^{2g}$. It is well-known that 
$\Diff^+(L)$ acts on $H^1(L,\Z)\cong \Z^{2g}$ as the 
integral $2g\times 2g$ symplectic group $Sp(2g,\Z)$. 
For each 
$0\neq \s \in H^1(\sg,\Z)$ there is a unique 
positive integer $m(\s)$ such that $\s =m(\s)\s '$ where 
$\s '\in H^1(\sg,\Z)$ is primitive. We call $m(\s )$ 
the {\em multiplicity} of $\s$. The multiplicity of 
$0\in H^1(\sg ,\Z)$ is defined to be 0.

Obviously,  $H^1(L,\Z)$ is trivial if $L=S^2$ is 
the 2--sphere. In this case $\la(S^2, \w', \w )=0$ is trivial. 

\begin{defn} 
{\rm 
Assume that the genus of $L$ is positive. 
Let $\s \in 
H^1(L, \Z)$ be the class representing
$\W_{\w',\w}$. Define 
\[ 
\la (L,\w',\w):=m(\s ) 
\] 
}
\end{defn}

\begin{lem} 
$\la (L,\w',\w)=\la (L,\w,\w')$. 
\end{lem} 

%
%
%

\begin{rem} 
{\rm 
When 
$L =T^2$ is of genus 1, then $TT^2$ and 
$T^*T^2$ are trivial $\R^2$-bundles over $T^2$. 
With trivializations of $TT^2$ and $T^*T^2$ fixed, 
the orientation-preserving bundle isomorphism 
$\W_\w:=TT^2\to T^*T^2$ then induces a map 
$T^2\to GL^+(2,\R)$ and hence an absolute invariant 
$\la (L,\w)\in \N\cup \{ 0\}$ is also defined.  
In particular we have $\la (T^2,\wc )=0$. 
This invariant $\la (T^2,\w)$ 
is actually the $\la$--invariant defined by 
Fintushel and Stern \cite{FS}. Our 
construction here can be thought as a relative extension 
of the Fintushel-Stern invariant to Lagrangian surfaces 
of any positive genus. 
} 
\end{rem}

\subsection{Realization of $\la (L,\w',\w )$} 

Consider the projection $\pi:T^*L\to L$, $(x,y)\to x$, 
and its differential 
\begin{equation}   
d\pi (x,y):T_{(x,y)}T^*L\to T_xL . 
\end{equation}   
For each $\rho:T^*L\to T^*L$ which is 
an orientation-preserving linear  
bundle automorphism over $L$, we define a 1-form 
$\la_{\rho}\in \W^1(T^*L)$ by 
\[ 
\la _{\rho}(x,y)=\rho(y)\circ d\pi (x,y):T_{(x,y)}T^*L\to \R 
\] 
and denote its negative differential $\w_{\rho}:=-d\la_{\rho}$. 
In particular, if $\rho=id$ then $\w_{\rho}=\wc$ the canonical 
symplectic 2--form on $T^*L$. Note that for each $\rho$ as 
defined above, 
$\w_{\rho}$ is symplectic and has the zero section $L$ as a 
Lagrangian surface.

Here is a way of constructing $\w$ from $\wc$ and a 
smooth map $\rho:=L\to S^1$.  
Fix a Riemannian metric $g_o$ on $L$. $g_o$ 
induces a Riemannian metric on $T^*L$. 
For any $\phi\in\text{Diff}^+(L)$ and any smooth map 
$\rho:L\to S^1$ consider the diffeomorphism 
\[ 
\Phi :T^*L \to T^*L 
\] 
\[ 
\Phi (x,y):=(\phi(x),e^{2\pi i\rho (\phi(x))}(\phi^{-1})^*y) , 
\quad p\in T^*_qL 
\] 
Here we use the polar coordinates for $T_x^*L\cong \C$. 
Clearly $\Phi(L)=L$. Let $\w:=\Phi^*\wc$. Then $L$ is also 
$\w$--Lagrangian. 
Then $\la(L,\w,\wc)$ is the multiplicity of the element 
in $H^1(L,\Z)$ corresponding to $\rho$. 
Note that $\w=\wc$ if $\rho$ is a constant map.

\subsection{The definition of $\la(L,L_0;\phi)$}

The construction of $\la (L,\w',\w)$ in Section \ref{local-1} 
can be applied straight forwardly 
to define an invariant $\la (L,L_0,\phi)$ for a pair of 
diffeomorphic Lagrangian surfaces $L_0$, $L=\phi(L_0)$ 
in a symplectic 4--manifold $(M,\w)$ and $\phi\in\Diff^+(M)$. 
We assume that $L_0$ is embedded and compact. 
Let  $\w_1:=(\phi^{-1})^*\w$. Then we define 
\[ 
\la (L,L_0,\phi):=\la (L,\w_1,\w) 
\] 
where $\la (L,\w_1,\w)$ is defined by using the fact 
that a tiny tubular neighborhood of $L$ is symplectomorphic to 
a tiny tubular neighborhood of the zero section of $T^*L$. 
This is really the {\em Lagrangian Neighborhood Theorem} due 
to Weinstein \cite{W}. We state the theorem here for future 
reference (see \cite{MS}). 

\begin{theo} \label{lag} 
Let $(M,\w)$ be a symplectic manifold and $L\subset M$ a compact 
Lagrangian submanifold. Then there exists a neighborhood 
$U\subset T^*L$ of the zero section, a neighborhood 
$V\subset M$ of $L$, and a diffeomorphism $\Phi:U\to V$ 
such that 
\[ 
\Phi^*\w=-d\lac, \quad \Phi|_L=id, 
\] 
where $\lac$ is the canonical 1--form on $T^*L$. 
\end{theo}

\subsection{Proof of Theorem \ref{lambda}} 

We start with the following lemma which is an easy consequence 
of the definition of $\la (L,L_0;\phi)$. 

\begin{lem} 
If $\phi \in \Diff^+_o(M;L_0\to L)$ is homotopic 
in $\Diff^+_o(M;L_0\to L)$ to 
some $\psi\in \Symp(M,\w)$ then 
$\la (L,L_0;\phi)=0$. 
\end{lem} 

The following Lemma will be frequently used in the proof of 
Theorem \ref{lambda}.

\begin{lem}[Lemma 3.14 of \cite{MS}]  \label{nbd} 
Let $M$ be a $2n$--dimensional smooth manifold and $L\subset M$ 
be a compact submanifold. Suppose that $\w_0,\w_1\in \W^2(M)$ 
are closed 2--forms such that at each point $q$ of $L$ the 
forms $\w_0$ and $\w_1$ are equal and nondegenerate on $T_qM$. 
Then there exist open neighborhoods $\cN_0$ and $\cN_1$ of $L$
and a diffeomorphism $\psi:\cN_0\to \cN_1$ such that 
\[ 
\psi |_L=id, \quad \psi^*\w_1=\w_0 . 
\] 
\end{lem}

Let us also recall Moser's argument on the isotopy of symplectic 
forms (\cite{Mo},\cite{MS}). Here we follow the presentation 
in \cite{MS}. For every family of symplectic forms $\w_t\in 
\W^2(M)$ with an exact derivative $\dfrac{d\w_t}{dt}=d\s_t$ 
there exists $\phi_t\in \Diff(M)$ such that
$\phi^*_t\w_t=\w_0$. $\phi_t$ can be chosen to be the time $t$ 
map of the flow of a time dependent vector field $X_t$, i.e., 
\[ 
\frac{d}{dt}\phi_t=X_t\circ \phi_t, \quad \phi_0=id, 
\]   
where $X_t$ satisfies the equation 
\[ 
d(\s_t+\io(X_t)\w_t)=0
\] 
In particular, since $\w_t$ is nondegenerate there exists a 
unique $X_t$ satisfying 
\[ 
-\s_t=\io(X_t)\w_t 
\] 

\begin{lem} \label{la=0} 
Assume that $\la (L,L_0;\phi)=0$ then $\phi$ can be smoothly 
isotoped to $\psi$ with $L$ fixed by the isotopy such that 
\[ 
(\psi^{-1})^*\w=w  \quad \text{ near } L_0
\] 
\end{lem} 

\begin{proof} 
Recall that 
\[ 
T_LM\cong T_L(T^*L)=TL\oplus T^*L
\] 
Since Since $\la (L,\w_1,w)=\la(L,L_0;\phi)=0$,  
\[ 
\rho:=\W_{\w_1,\w}: T^*L\to T^*L  
\] 
is homotopic to to the identity map. 
Fix a trivialization of $TL\otimes T^*L$ then 
$\rho$ becomes a map $L\to GL^+(2,\R)$ and is 
represented by $A:=(a_{ij})\in GL(2,\R)$ 
which depends smoothly on $x\in L$. 
We may assume that $\w=-d\lac$ near $L$. 
Since $L$ is $\w_1$--Lagrangian, $\w_1$ is exact near $L$. 
Let $f\in C^{\infty}(M)$ be a smooth function 
supported in a small tubular neighborhood of $L\subset M$ 
and $f=1$ near $L$. Define 
\[ 
\tau :=\w_1+d(f\rho(y)\circ d\pi) 
\] 
$\tau\in\W^2(M)$ is a closed 2--form which is exact near $L$ 
and $\tau|_{TL}=0$. Near $L$ we have $\tau=d\h$ for some 
1-form $\h$ defined near $L$. Let $i:L\hookrightarrow M$ 
denote the inclusion of $L$. Then $i^*\h$ is a closed 1-form 
on $L$. Then near $L$ we have 
\[ 
\w=-d\lac, \quad \w_1=-d(\rho(y)\circ d\pi )+d\h  
\]

Let $\rho_t:T^*L\to T^*L$, $t\in [0,1]$, be a smooth homotopy 
of linear automorphisms between $\rho_0:=\rho$ and
$\rho_1:=id$.  By using Theorem \ref{lag} we 
can    
consider a family of symplectic 2-forms defined near $L$: 
\[ 
\w_t:=-d(\rho_t(y)\circ d\pi)+td\h . 
\] 
Consider the vector field $Y_t$ defined by 
\[ 
-\frac{d\rho_t}{dt}\circ d\pi +\h =\io(Y_t)\w_t 
\] 
Note that $Y_t=0$ on $L$ for all $t$. 

Now consider a smooth 
time dependent vector field $X_t$ supported in a tubular 
neighborhood of $L$ such that on $X_t=Y_t$ on $L$.   
Define $\phi'_t$ to be the time $t$ map of the flow of $X_t$. 
$\phi'_t$ fixes $L$ pointwise for all $t$, $\phi'_0=id$. 
We obtain an  isotopy $\phi'_t\circ\phi$ between 
$\phi$ and $\phi'_1\circ\phi$ such that $\psi':=\phi'_1\circ\phi$ 
satisfies 
\[ 
(\psi')^*\w =\w \quad \text{ on } L_0 
\] 
By Lemma \ref{nbd} we can further isotope $\psi'$, with $L$ 
fixed, to a diffeomorphism $\psi$ such that $\psi(L_0)=L$ and 
$\psi^*\w=\w$ near $L_0$. 
\end{proof}

\begin{lem} \label{symp} 
Let $\phi$ be as in Lemma \ref{la=0}. Assume that 
$\la (L,L_0;\phi)=0$ and  
$\phi\in\Diff^+_0(M)$, then there exists  
$\Phi_t\in \Diff^+_0(M;L_0\to L)$, $0\leq t\leq 1$, such that 
\[ 
\Phi_0=\phi, \quad (\Phi_1^{-1})^*\w=\w . 
\] 
That is, $L$ is symplectomorphic to $L_0$. 
\end{lem}

\begin{proof} 
Let $\phi_t\in\Diff^+_0(M)$, $0\leq t\leq 1$, be a smooth 
homotopy between $\phi_0:=id$ and $\phi_1:=\phi$. Denote 
$\w_t:=(\phi^{-1}_t)^*\w$. $w_0=w_1$ near $L$.  
Since $\la (L,L_0;\phi)=0$, by Lemma \ref{la=0} 
we may assume that 
$\w_t=\w_0$ and $\w_{1-t}=w_1$ for all $0\leq t\leq \de_o$ 
for some constant $\de_o>0$.  $L$ 
need not be $\w_t$--Lagrangian for $t\in (\de_o,1-\de_o)$.

Write 
$w_t=\w +d\alpha_t$ with $\alpha_t\in \W^1(M)$, 
$\alpha_0=0$. Since $w_1=w_0$ near $L$, $d\al_1=0$ 
near $L$.  
Let $i:L\hookrightarrow M$ denote the embedding of $L$,  
and $\s_t:=i^*\al_t$. Then $\s_0=0$, 
$\s_1$ is a closed 1-form on $L$. \\ 
 
\noindent 
{\bf Claim:} There exists a closed 1-form $\be\in \W^1(M)$ 
such that $i^*\be=\s_1$. 

\begin{proof}[Proof of the Claim] 
\quad \\ 
Consider the long exact sequence of cohomology groups 
\begin{equation}  \label{ex-seq} 
\cdots \to H^1(M)\overset{i^*}{\to} 
 H^1(L)\overset{\de^*}{\to} H^2(M,L) \to \cdots 
\end{equation}  
Since $\de^*[\s_1]=[d\al_1]=0$, $[\s_1]=\io^*c$ for some 
$c\in H^1(M)$. Hence there exists a closed 1-form 
$\be'\in \W^1(M)$, $[\be']=c$, such that $i^*\be'=\s_1+df$ 
for some $f\in C^{\infty}(L)$. Extend $f$ to be a smooth 
function ( also denoted by $f$) on $M$. Define 
$\be:=\be'-df$. Then $\be\in \W^1(M)$ is closed and 
$i^*\be=\s_1$. 
\end{proof} 

Subtract $t\be$ from $\al_t$ and still call the resulting 
1--form $\al_t$. Then
\[ 
\w_t=\w_0+d\al_t, \quad \al_0=0, \quad 
\s_0=i^*\al_0=0, \quad \s_1=i^*\al_1=0. 
\] 

Recall the projection $\pi:T^*L\to L$. Let $h$ be a 
smooth function on $T^*L$ that is compactly supported in 
a tiny tubular neighborhood of the zero section $L$ and 
satisfies $h=1$ near $L$. By Lemma \ref{nbd} the 1--form 
$\g_t:=h\pi^*\s_t$ is a 1--form supported in a tubular 
neighborhood of $L\subset M$ such that $i^*\g_t=\s_t$. 
For each $t$ consider the time independent vector field 
$X_{t}$ defined by 
\[ 
\g_t=\io(X_{t})\w_t 
\] 
Let $\psi_{t,s}$ denotes the time $s$ map of the flow of 
$X_t$. $\psi_{0,s}=id=\psi_{1,s}$ for all $s$. 
Denote $\w_{t,s}:=(\psi_{t,s}^{-1})^*\w_t$. 
Then $i^*\w_{t,1}=0$ for all $t$, i.e., $L$ is 
$\w_{t,1}$--Lagrangian. Write $\w'_t:=\w_{t,1}$. On $L$, 
as discussed in Lemma \ref{la=0},  
\begin{align*}  
\w'_t&=-d(\rho_t(y)\circ d\pi)+d\h_t, \quad w'_0=w'_1 \\ 
i^*\h_t &\in \W^1(L) \text{ is closed}, \quad 
\rho_0(y)=\rho_1(y)  
\end{align*}  

$\phi'_t:=\psi_{t,1}\circ\phi_t$, 
$t\in [0,1]$, is a homotopy between $id=\phi'_0$ and 
$\phi=\phi'_1$. $(\phi'_t)^*\w'_t=w_0$,  
$\w'_0=\w$, $\w_1'=\w_1$. $\phi'_1(L_0)=L$, 
$L$ is $\w'_t$--Lagrangian $\forall t\in [0,1]$.  
Note that $\la (L,L_0;\phi'_t)=0$ for all $t$. 

Again, by applying  Lemma \ref{la=0} to each of the triples  
$(L,L_0,\phi'_t)$, $t\in [0,1]$, the homotopy 
$\phi'_t$ is perturbed, with $L$ fixed, to a new homotopy 
$\phi''_t$ between $\phi''_0=id$ and $\phi''_1=\phi$  
such that  
\[ 
\w''_t:=((\phi''_t)^{-1})^*\w \ \ \text{ equals } \w 
\text{ near } L 
\] 

Hence $\w''_t=\w+d\eta_t$ with $\eta_0=0$ and 
$\eta_t\in \W^1(M)$ closed near $L$. 
Note that $\w''_1=\w_1$. By applying the 
long exact sequence (\ref{ex-seq}) we may assume that 
\[ 
\eta_t=0 \quad \text{ near } L.  
\] 
Define the vector field $Y_t$: 
\[ 
\frac{d\eta_t}{dt}=\io(Y_t)\w''_t 
\] 
and let $\psi_t$ be the time $t$ map of the flow of $Y_t$. 
$\psi_0=id$, $\psi _t(L)=L$ for $t\in [0,1]$, 
and $\psi_1^*\w=\w''_1=\w_1$. 
Note that 
\[ 
(\psi_1\circ\phi)^*\w=\phi^*\psi^*_1\w=\phi^*\w_1=w 
\] 
Then $\Phi_t:=\psi_t\circ \phi\in\Diff^+_o(M;L_0\to L)$ is an  
homotopy between $\Phi_0=\phi$ and 
$\Phi_1=\psi_1\circ\phi\in \text{Symp}(M,\w)$. 
This completes the proof of Lemma \ref{symp} and 
hence the proof of Theorem \ref{lambda}. 
\end{proof}

%
%
%

\section{A smooth isotopy invariant} \label{second}

\subsection{Local construction}

Let $\tau$ be an area form on a closed oriented Riemann surface 
$L$. For $k_1,k_2\in\R$ with $k_1k_2<0$,  
$k_1\tau\oplus k_2\tau$ is a symplectic 2--form on 
$L\times L$. Let 
$\Delta:=\{ (x,x)\mid x\in L\}$ 
denote the diagonal of $L\times L$. $\Delta$ is 
$(k_1\tau\oplus k_2\tau)$--Lagrangian iff $k_1+k_2=0$, 
$(k_1\tau\oplus k_2\tau)$--symplectic iff $k_1+k_2>0$. 
By allpying Lemma \ref{lag}, 
there exists a diffeomorphism $\phi:\cN (L,T^*L)\to 
\cN(\Delta,L\times L)$ from a tubular neighborhood 
$\cN (L,T^*L)$ of the zero section of $T^*L$ to 
a tubular neighborhood $\cN(\Delta,L\times L)$ 
of $\Delta \subset L\times L$ such that $\phi^*
(k_1\tau\oplus (-k_1)\tau )=\wc$. 

We can get a second symplectic 
2-form on $T^*L$ as follows: By scaling the fibers of $T^*L$ 
we get a fiber-preserving diffeomorphism $\psi:T^*L \to 
\cN (L,T^*L)$. Fix constants $k_1,k_2$ with 
$k_1k_2<0$ and $k_1+k_2>0$, 
then the pullback $\w^+:=(\phi\circ\psi)^*
(k_1\tau\oplus k_2\tau)$ is a symplectic 2-form. 

Fix a Riemannian metric $g$ on $T^*L$ for example, let 
$g$ be one that is induced by a Riemannian metric on $L$. 
Then we get an $S^2=SO(4)/U(2)$--bundle $\cE$ over $T^*L$, 
whose sections are in one-one correspondence with 
{\em $g$--skew-adjoint} ($g$--skad)  
almost complex structures on $T^*L$. 
Recall from \cite{MS} (p.60--62)
that for each symplectic 2--form 
$\w$ on a Riemannian manifold $(M,h)$ there associates a 
unique $\w$--compatible almost complex structure that is 
also $h$-skad. Now that we have two linearly 
independent symplectic 2--forms $\wc$ and $\w^+$ on 
$(T^*L,g)$, each of $\wc$ and $w^+$ associates a unique 
$g$--skad almost complex structure denoted by $J_g$, $J_{g+}$ 
respectively. $J_g$ and $J_{g+}$, viewed as two sections of 
$\cE$, never intersect, hence the triple sections 
$(J_g, J_{g+}, J_{g+}\circ J_g)$ defines  a trivialization 
\[ 
\cE\cong S^2\times T^*L. 
\] 

Let $L'\subset T^*L$ be an embedded $\wc$--Lagrangian 
surface. Assume that there exists $\phi\in 
\Diff ^+_o(T^*L)$ 
such that $\phi(L)=L'$ and $\phi$ is $\wc$--symplectic near 
$L$, i.e., $\phi^*\wc=\wc$ in a tubular neighborhood of $L$. 
Take a smooth path $\phi_t\in \Diff ^+_o(T^*L)$ with 
$\phi_0=id$ and $\phi_1=\phi$. Define $\w_t:=(\phi^{-1})^*\wc$ 
and $L_t:=\phi_t(L)$. We have that $L_t$ is 
$\w_t$--Lagrangian, $L_0=L$, $L_1=L'$, and $\w_0=\wc$. 
Note in particular $w_1=\wc$ near $L'$. 

Let $J_t$ denote the $\w_t$-compatible 
$g$--skad almost complex structure. Think of the family $L_t$, 
$0\leq t\leq 1$, as the image of a map $\eta :L\times [0,1]\to 
T^*L$. $\eta^*\cE$ is trivial $S^2$--bundle over $L\times 
[0,1]$ with the trivialization induced by $(J_g, J_{g+},
J_{g+}\circ J_g)$. Now the union 
$J(\eta):=\cup _t\eta^*(J_t|_{L})$ is a 
an element of 
\[ \Gamma(\eta^*\cE,J_g)=\{ J\in\Gamma(\eta^*\cE)\mid 
J=\eta^*J_g \ \text{ on } \ L\times \{ 0,1\}\} ,  
\] 
the set of all $C^{\infty}$ sections of $\eta^*\cE$ which 
are equal to $\eta^*J_g$ over $L_0\cup L_1$. 

Clearly, a necessary condition for 
the path $L_t$ homotopic to a $\wc$--Lagrangian  isotopy  
relative to $(L_0,L_1)$ is that $J(\eta)$ is homotopic to 
the "constant" section $\eta^*J_g$ in 
$\Gamma(\eta^*\cE,J_g)$, i.e., 
if $[J]$ is in the connected component of 
$\pi_0(\Gamma(\eta^*\cE,J_g))$ 
containing $\eta^*J_g$.

We can think of $J_g$ as the north pole $*$ of the 2--sphere of
almost complex  structures of the given Riemannian metric $g$.
The set 
$(S^2)^L$ of continuous maps 
from $L$ to $S^2$ has a special constant map 
denoted by  $*_L$ which sends 
$L$ to the point $*$. Then 
\[ 
\Gamma(\eta^*\cE,J_g)\hookrightarrow (S^2)^L_o
\] 
where $\cM:=(S^2)_o^L$ is the connected component of 
$(S^2)^L$ containing the constant map $*_L$.  
Let $\W\cM$ denote the loop space of 
the pointed space $(\cM,*^L)$. 

Given $\{ f_t\}\in \W\cM$, i.e., 
$f_t:L\to S^2$, $t\in [0,1]$, is a smooth path in 
$\cM$ with $f_0=*_L=f_1$. We can associate a map 
\[ 
\hat{F}:L\times [0,1] , \quad \hat{F}(\cdot, t)=f_t 
\] 
Since $f_0=*_L=f_1$, $\hat{F}=F\circ u$ where 
\[ 
u:=L\times [0,1] \to X =L\times [0,1]/L\times \{ 0,1\} 
\] 
$X$ is the double cone of 
$L$ with the two vertices identified, and $u$ is 
the crushing map that sends $L\times \{ 0, 1\}$ to the 
vertex of $X$, and $u|_{L\times (0,1)}$ is a homeomorphism. 
Then each $\{ f_t\}\in \W\cM$ induces a map 
\[ 
F:X\to S^2 
\] 
It is easy to see that if $\{ f_t\},\{ f'_t\}\in \W\cM$ 
are homotopic, then the corresponding maps 
$F,F':X\to S^2$ are homotopic as continuous 
maps from $X$ to $S^2$. 

Think of $J=\{ J_t\}$ as a loop $\{ f_t\}$ in $\cM$ and hence 
a map from $X\to S^2$, we arrive at the following  

\begin{defn} 
{\rm 
$n(J)$ is defined to be the corresponding homotopy class in 
$[X,S^2]$ that is represented by $F$. 
}
\end{defn}

Let $U:=L\times (0,1)$ and $V:=X\setminus L\times\{ 
\frac{1}{2}\}$. Then $X=U\cup V$. 
By considering the Mayer-Vietoris sequence of $X=U\cup V$ 
we get a 
long exact sequence of homology groups (with $\Z$-coefficients) 
\[ 
\cdots \to H_i(U\cap V)\to H_i(U)\oplus H_i(V)\to 
H_i(X)\to H_{i-1}(U\cap V)\to \cdots 
\] 
and get 
\[ 
H_1(X)\cong \Z, \quad H_2(X)\cong \Z^{2g}, \quad H_3(X)\cong 
\Z 
\] 
hence the cohomology groups 
\[ 
H^1(X)\cong \Z, \quad H^2(X)\cong \Z^{2g}, \quad H^3(X)\cong 
\Z 
\] 

The following lemma is due to Kuperberg \cite{K}: 

\begin{lem} 
The homotopy class of a map $F$ from $X$ to $S^2$ is described 
by $c:=F^*\s$, where $\s\in H^2(S^2,\Z)$ is the positive 
generator, and a Hopf degree $d\in \Z_n$ where 
$n$ is the maximal divisor of $c$. 
\end{lem} 

\begin{cor}
If $L=S^2$, then $[X,S^2]\cong \Z$. 
\end{cor} 

\begin{exam} 
{\rm 
Think of $S^3$ as the reduced suspension of $S^2$, 
then the Hopf map $H:S^3\to S^2$ gives a loop $f_t$  
in $(S^2)^{S^2}_o$ with $f_0=*_{S^2}=f_1$ and 
$f_t(p)=*$ for some $p\in S^2$, $\forall t\in [0,1]$. 
The loop $f_t$ represents a nontrivial element of 
$\pi_1((S^2)^{S^2})$. More generally, for any 
degree $d$ map $\Phi_d:S^3\to S^3$, the map $H\circ \Phi_d$ 
induces an element of $\pi_1((S^2)^{S^2})$, and two 
such elements are distinct if $d_1\neq d_2$. 
} 
\end{exam} 

\begin{rem} 
{\rm 
When the genus of Lagrangian surfaces are positive, 
$n(L_1,L_0;[L_t])$ indicates the potential existence  
of symplectomorphisms 
which are smoothly but not symplectically isotopic to the 
identity map, and which are {\em not} generalized Dehn twists. 
Are there such symplectomorphisms? It will be very 
interesting if one can construct such an example or 
disprove the existence. 
} 
\end{rem}

Note that a priori $n(J)\neq 0$ 
does not guarantee that $L_1$ is not $\wc$--Lagrangian 
isotopic to $L_0$, since in general 
$n(J)$ depends on the choice of 
the path $L_t$.

\subsection{Definition of $n(L_1,L_0;[L_t])$}

Let $L_t=\phi_t(L_0)$, $\phi_t\in \Diff^+_o(M)$, 
$t\in [0,1]$, be a smooth isotopy of embedded surfaces 
in a symplectic 4--manifold $(M,\w)$ with $L_0,L_1$ 
being $\w$--Lagrangian and 
$\phi_0=id$, $\phi_1^*\w=\w$ near $L_0$ (so $L_0,L_1$ 
are symplectomorphic).  
Fix an $w$--compatible almost complex 
structure $J_\w$ and let $g_\w=\omega\circ(Id\times J)$ 
be the induced Riemannian metric on $M$.  
Let $\cE$ denote the associated $S^2$--bundle over $M$. 
Fix a symplectic 2--form $\w_+$ which is defined near 
$L_0$ such that $L_0$ is $\w^+$-symplectic. 
Define $\w_t:=(\phi^{-1}_t)^*\w$, $\w^+_t:=(\phi^{-1}_t)^*\w^+$. 
$L_t$ is $\w_t$--Lagrangian and $\w^+_t$--symplectic. 
Let $J_t$ (resp. $J^+_t$) denotes the 
$\w_t$--compatible (resp. $w^+_t$--compatible) $g_\w$-skad almost
complex structure on $T_{L_t}M$. Then the triple 
$(J_t,J^+_t,J^+_t\circ J_t)$ trivializes the $S^2$-bundle  
$\cE|_{L_t}\cong S^2\times L_t$.  
Note that $J_t|_{L_t}=J_\w$ for $t=0,1$. 
With respect to the $(J_t,J^+_t,J^+_t\circ J_t)$ 
trivialization, 
the section $J_\w|_{L_t}$ becomes a map 
$L\times [0,1]\to S^2$ with $n(J_w,[\phi_t])$ 
is 0 iff $J_\w$ is homotopic to $J_t$ as a section of 
$\cE|_{L_t}$ with boundary fixed.

\begin{defn} 
{\rm 
$n(L_1,L_0;[L_t])=n(L_1,L_0;[\phi_t]):=n(J_w,[\phi_t])$. 
}
\end{defn} 

The following lemma is a easy consequence of the definition 
of $n(L_1,L_0;[\phi_t])$. 

\begin{lem} 
If the path $L_t:=\phi_t(L_0)$ ($\phi_0=id$) is 
homotopic to 
a Lagrangian isotopy with boundary $L_0\cup L_1$ preserved, 
then $n(L_1,L_0;[\phi_t])=0$. 
\end{lem}

\begin{ques} 
{\rm 
Is it true that $n(L_1,L_0;[L_t])=0$ implies that 
the path $L_t:=\phi_t(L_0)$ ($\phi_0=id$) is 
homotopic to 
a Lagrangian isotopy with boundary $L_0\cup L_1$ preserved? 
} 
\end{ques}

Any answer to the above question will help us better 
understand the isotopy problem of Lagrangian surfaces. 
Also, with Theorem \ref{lambda} proved, it will be very 
important to construct or find examples Lagrangian surfaces 
(of positive genus) 
which are diffeomorphic but not symplectomorphic. 
Any such example will shed some new light on our 
understanding of symplectic topology. 
Moreover, by applying the philosophy behind the construction 
of $\la (L,L_0;\phi)$ and $n(L_1,L_0;[L_t])$ to 
the contact case, we can also define similar invariants 
for Legendrian knots in contact 3--manifolds, and use 
them to explore Legendrian isotopy problems \cite{Y}. 
We will come back to these topics later.

\end{document}